\documentclass[a4paper, oneside, 10pt]{article}
\RequirePackage[OT1]{fontenc}
\RequirePackage{amsthm,amsmath}
\RequirePackage[numbers]{natbib}
\usepackage{graphics}
\usepackage{latexsym}
\RequirePackage[colorlinks,citecolor=blue,urlcolor=blue]{hyperref}
\usepackage{amsfonts}
\usepackage{amssymb}
\usepackage[latin1]{inputenc}
\usepackage[a4paper, margin=3.4cm]{geometry}
\usepackage{verbatim}

\theoremstyle{plain}

\newtheorem{prop}{Proposition}[]
\newtheorem{coro}{Corollary}[]
\newtheorem{rema}{Remark}[]

\theoremstyle{definition}

\begin{document}
	
	{\linespread{1}
		\title{{\bf{\LARGE{On two overlooked stick-breaking constructions of the normalized inverse Gaussian process}}}\footnote{{\it AMS (2000) subject classification}. Primary: 60G58. Secondary: 60G09.}}
		\author{\textsc {Annalisa Cerquetti}\footnote{Corresponding author. e-mail: {\tt annalisa.cerquetti@gmail.com}}\\
			\it{\small Independent researcher. Rome, Italy}}
		\date{\today}
		\maketitle{}
	}
	
	\begin{abstract}
		We shed light on two alternative stick-breaking constructions of the normalized inverse Gaussian (NIG) random discrete distribution
		which appear to have been overlooked so far in the Bayesian nonparametric setting.
		The first is derived from a result in Aldous and Pitman (1998) for the conditional Brownian excursion partition, mixing over the local time at zero up to time one. 
		The second arises as a particular case of a result in James (2013) for priors obtained by a random spatial and temporal change of the normalized generalized Gamma subordinator. Both constructions are in terms of straightforward transformations of standard random variables and can be easily generalized to provide the stick-breaking construction of any element,  respectively, in a) the family of mixed Poisson-Kingman models driven by the $1/2$ stable L\'evy measure and b) the family of Poisson-Gamma processes driven by the Inverse Gaussian subordinator.  
		
	\end{abstract}

	\section{Introduction}
	
	The normalized inverse Gaussian (NIG) random discrete distribution was originally introduced in Pitman (2003, Sect. 5.4) as an instance of a tractable {\it mixed} Poisson-Kingman model PK($\rho, \gamma$) driven by the stable $(1/2)$ L\'evy density $\rho_{1/2}(s)=\frac{1}{2\sqrt{\pi}}s^{-3/2}$,
	with mixing distribution $\gamma$ the density of an {\it inverse Gaussian} random variable of parameters $(1/v, 1)$ for $v>0$, namely
	\begin{equation}
		\label{IGdensity}\gamma_{IG, v}(y)=\sqrt{\frac{1}{2\pi}} y^{-3/2} \exp \left\{-\frac{1}{2}\left(v^2y +y^{-1}-2v\right)\right\}dy.
	\end{equation}
	Here {\it mixed} Poisson-Kingman models are distributions for {\it ranked} random discrete distributions, $(P_i^{\downarrow}):=(J_i^{\downarrow}/T)$,  derived from an inhomogenous Poisson point process of random lengths $(J_1^{\downarrow} \geq  J_2^{\downarrow} \geq \dots \geq 0 )$ by normalizing by the total length $T:=\sum_{i=1}^\infty J_i^{\downarrow}$, conditioning on $T=t$ and mixing the distribution of $(P_i^{\downarrow}|t)$ over $t$ with respect to a general density $\gamma$.
	In Bayesian nonparametrics  the normalized inverse Gaussian prior has been first investigated in Lijoi {\it et al}. (2005) as a tractable alternative to the Dirichlet prior (Ferguson, 1973) in hierarchical mixtures modeling. In the same setting, relying on the general theory in Perman {\it et al}. (1992) for the size-biased sampling of the jumps of a subordinator,  Favaro {\it et al}. (2012) derived a stick-breaking construction of the size-biased permutation $(\tilde{P}_j)_{j \geq 1}$ of its ranked atoms showing that, for $G_1$ a {\it generalized inverse Gaussian} (GIG) r.v. of parameters $\left(v, 1,-\frac{1}{2}\right)$, $v >0$ and $Z_1$ a {\it L\'evy} $(0,1)$ r.v., then 
	\begin{equation}
		\label{struFAV}
		\tilde{P}_{1, NIG}\stackrel{d}{=}\frac{G_1}{G_1+Z_1} 
	\end{equation}
	with density 
	\begin{equation}
		\label{p1fav}
		f_{\tilde{P}_{1, NIG}}(p)=\frac{v^{1/4}}{(2p\pi)^{1/2} (1-p) K_{-1/2}(v^{1/2})}K_{-1}\left\{\left(\frac{v}{1-p}\right)^{1/2}\right\},
		\\
	\end{equation}
	for $0 < p < 1$, $K_{-1}=K_{1}$ and $K_1$ the usual modified Bessel function of the {\it second} kind 
	$$
	K_1(u)=\frac{1}{2u}\int_0^\infty \frac{1}{x^2} \exp\left\{ -\frac{1}{2} \left(u^2x +\frac{1}{x}\right)\right\}dx.
	$$
	Additionally, for $j \geq 2$, they proved that $\tilde{P}_{j, NIG}\stackrel{d}{=}V_j\prod_{i=1}^{j-1} (1-V_i)$, for 
	\begin{eqnarray}
		\label{vjFAV}
		\left(V_j \mid V_1, \ldots, V_{j-1}\right)\stackrel{}{=}\frac{G_j}{G_j+Z_j}  
	\end{eqnarray}
	where
	\begin{equation}
		\label{vjFAV2}
		G_j \sim \operatorname{GIG}\left(\frac{v}{\prod_{i=1}^{j-1}\left(1-V_i\right)}, 1,-\frac{j}{2}\right)
	\end{equation}
	and the $(Z_j)_{j \geq 1}$ are independent L\'evy $(0,1)$  r.v.'s, independent from the $(G_j)_{j\geq 1}$. 
	Actually, in Pitman (2003, Sect. 5.4),
	the {\it structural distribution} - i.e. the law of  the first size-biased pick $\tilde{P}_{1, NIG}$ from $P \sim PK(\rho_{1/2}, \gamma_{IG, v})$, - 
	has been shown to admit the simpler density 
	\begin{equation}
		\label{besselp1}
		f_{\tilde{P}_{1, NIG}}(p) = \frac{v \exp(v)}{\pi \sqrt{p}(1-p)}K_1\left(\frac{v}{\sqrt{1-p}}\right),
	\end{equation}
	for $v>0$, and $0<p<1$.  And indeed (\ref{besselp1}) can be derived from (\ref{p1fav}) recalling that, by the theory of Bessel functions, $K_{1/2}(x)=\exp(-x)\sqrt{\frac{\pi}{2}}\frac{1}{\sqrt{x}}$.
	
	Note also that, by an easy manipulation, the density of a Beta ($1/2, 1/2$) r.v. can be isolated in (\ref{besselp1})
	\begin{equation}
		\label{pit_bes}
		f_{\tilde{P}_{1, NIG}}(p) = \frac{1}{\pi \sqrt{p(1-p)}}\exp(v)\frac{v}{\sqrt{1-p}}K_1\left(\frac{v}{\sqrt{1-p}}\right),
	\end{equation}
	and, expliciting the Bessel function, (\ref{pit_bes}) yields
	\begin{equation}
		\label{pit_ig}
		f_{\tilde{P}_{1, NIG}}(p) = \frac{1}{\pi \sqrt{p(1-p)}}\exp(v)\frac{1}{2}\int_0^\infty \frac{1}{y^2}\exp\left\{ -\frac{1}{2} \left(\frac{v^2 y}{1-p} +\frac{1}{y}\right)\right\}dy.
	\end{equation}
	
	Here we shed light on two further stick-breaking representations of the size-biased permutation of the ranked atoms of the normalized Inverse Gaussian random discrete distribution that - to the best of our knowledge - have remained unexplored so far in the Bayesian nonparametric literature. Their main feature is to provide characterizations in terms of straightforward transformations of standard random variables, therefore offering convenient alternatives to address computational issues arising with the implementation of simulation algorithms for posterior inference, both when dealing with hierarchical mixture modeling and species sampling problems (see e.g. Canale {\it et al}., 2022, Favaro {\it et al.}, 2014 and Favaro and Walker, 2013). Additionally their particularly simple structure allows easy generalizations to larger families of random discrete distributions, thus potentially enlarging the set of models that deserve to be investigated as tractable priors in Bayesian nonparametrics.  
	
	\section{Results}
	
	\subsection{A construction derived from Aldous and Pitman (1998) model}
	
	The first construction we are going to present is derived from a result for the Brownian excursion partition conditioned on the local time at zero up to time one (Aldous and Pitman, 1998, Th. 4, Corollary 5), diffusely discussed in Pitman (2003, Sect. 7 and 8), that here we summarize.

	Let  $\Pi$ be the {\it Brownian excursion partition}, that is the random partition of $N$ generated by uniform random sampling of points from the interval [0,1], partitioned by the excursion intervals of a standard Brownian motion $B$. Let $\{L_t, t>0\}$ be the local time at 0 of $B$ and recall that its inverse $\{\tau_\ell, \ell \geq 0\}$ is a $(1/2)$ stable subordinator, hence, by the usual normalization of the local time as occupation density relative to the Lebesgue measure, $\tau_1\stackrel{d}{=}2G$ for $G \sim  1/2-stable$ and $L_1 \stackrel{d}{=}|B_1|$. Then, by a result in Pitman and Yor (1992)  $\Pi$ is a PK($\rho_{1/2}$) partition with $1/2$-diversity $\sqrt{2}L_1$ and the sequence $(\tilde{P}_j)$ of class frequencies of $\Pi$ is the sequence of lengths of excursion intervals in order they are discovered by the sample points. Conditioning on $L_1=\lambda$, that is the local time of $B$ at $0$ up to time $1$,  for $\lambda \geq 0$ then 
	\begin{equation}
		\label{equiPK}
		\Pi(\lambda)\stackrel{d}{=}(\Pi|L_1=\lambda) \mbox{  is a  }PK(\rho_{1/2}| 1/2 \lambda^{-2})
	\end{equation} 
	and this is the common distribution of ranked lenghts of excursions of the Brownian bridge over $[0,1]$ given $L_1=\lambda$. Additionally, for $\tilde{P}_j(\lambda)$ the random relative frequency of the $j$-th class of $\Pi(\lambda)$ discovered in the process of random sampling, then $(\tilde{P}_j(\lambda), j=1, 2, \dots)$ are distributed as the lengths of excursions of $B$ over $[0,1]$ given $L_1=\lambda$ as discovered by a process of {\it length-biased sampling}.

	Now, again from Pitman (2003, eq. 94), the law of $\tilde{P}_1(\lambda)$, the structural distribution of $\Pi(\lambda)$, is given by  
	\begin{equation}
		\label{p1la}
		Pr(\tilde{P}_1(\lambda) \in dp) = \frac{\lambda}{\sqrt{2\pi} }p^{-1/2}(1-p)^{-3/2} \exp\left(-\frac{\lambda^2}{2}\frac{p}{1-p})\right)dp
	\end{equation}
	for $0 <p<1$ and $\lambda>0$, or equivalently 
	$
	P(\tilde{P}_1(\lambda)\leq p)= 2\Phi\left(\lambda \sqrt{\frac{p}{1-p}} \right)-1
	$
	for $\Phi$ the standard Gaussian distribution function, hence the following equality in distribution holds
	\begin{equation}
		\label{ratio}
		\tilde{P}_1(\lambda)\stackrel{d}{=}\frac{X_1}{\lambda^2 +X_1}\nonumber
	\end{equation}
	for $X_1 \sim \chi_1^2$ a chi-squared r.v. with one degree of freedom.
	It is easy to see that, by the previous discussion, a constructive definition of the structural distribution of the normalized inverse Gaussian ranked random discrete distribution can be straightforwardly derived as follows.
	
	{\prop 
		Let $X_1 \sim \chi_1^2$ and $T$ be an inverse Gaussian random variable of parameters $(1/v, 1)$, independent from $X_1$, with density (\ref{IGdensity}), then
		\begin{equation}
			\label{struPIT}
			\tilde{P}_{1, IG}\stackrel{d}{=}\frac{X_1} {T^{-1} + X_1}.
	\end{equation}}
	
	\proof {Recall that the inverse Gaussian density is obtained by exponentially tilting the density of $T=2G$ for $G$ a Stable ($1/2$) r.v. (Pitman, 2003, eq. 53), that is $T$ is a L\'evy $(0,1)$ r.v. with density
		$f_{T, 1/2}(t)=\frac{1}{\sqrt{2\pi}}t^{-3/2}e^{-\frac{1}{2t}}.$ 
		By (\ref{equiPK}) and the previous discussion,
		\begin{equation}
			\label{analogy}
			PK\left(\rho_{\frac12}\,\middle\vert\, \frac{1}{2\lambda^2}\right)\stackrel{d}{=}PK\left(\rho_{\frac12}\,\middle\vert\, \frac{t}{2} \right)\nonumber
		\end{equation}
		and, by the change of variable $1/t=\lambda^2$, (\ref{p1la}) yields 
		\begin{equation}
			\label{p_lambda}
			Pr(\tilde{P}_{1, IG}(t) \in dp) = \frac{1}{\sqrt{2\pi} }p^{-1/2}(1-p)^{-3/2} t^{-1/2}\exp\left(-\frac{1}{2t}\frac{p}{1-p}\right)dp.
		\end{equation}
		Mixing (\ref{p_lambda}) over $t$ with (\ref{IGdensity}) then, for $b=v^2$
		$$
		f_{\tilde{P}_{1, IG}}(p)=\frac{\exp\{\sqrt{b}\}}{\sqrt{2\pi}p^{1/2}(1-p)^{3/2} }\int_0^\infty t^{-2}\exp\left(-\frac{1}{2t}\frac{p}{1-p}\right) {\frac{1}{\sqrt{2\pi}}} \exp \left\{-\frac{1}{2}\left(bt +t^{-1}\right)\right\}dt
		$$
		$$
		=\frac{\exp\{\sqrt{b}\}}{{2\pi}}p^{-1/2}(1-p)^{-3/2}\int_0^\infty t^{-2}\exp\left(-\frac{1}{2}\left( bt + \frac{1}{t(1-p)}\right)\right)dt
		$$
		and, by the change of variable $y=t(1-p)$,  $dt=(1-p)^{-1}dy$, then
		$$
		=\frac{\exp\{\sqrt{b}\}}{{2\pi} }p^{-1/2}(1-p)^{-1/2}\int_0^\infty y^{-2}\exp\left(-\frac{1}{2}\left( \frac{by}{1-p} + \frac{1}{y}\right)\right)dy
		$$
		which agrees with (\ref{pit_ig}). \qed
	}

	The result in Proposition 1 can be readily extended to the construction of whole size-biased permutation $(\tilde{P}_j)_{j\geq 2}$ of the ranked atoms of the NIG process, again relying on a conditional result in Aldous and Pitman (1998, Corollary 5) recalled in  Pitman (2003, Proposition 14).

	\prop  Let $(X_i)_{i \geq 1}$ be a sequence of iid $\chi_1^2$ r.v.s,  $S_j:=\sum_{i=1}^j X_i$, and $T$ an Inverse Gaussian $(1/v, 1)$ r.v. independent from $(X_i)_{i \geq 1}$, then the size-biased permutation of the ranked atoms of the normalized Inverse Gaussian random discrete distribution admits the following  stick-breaking construction for $j\geq 2$ 
	\begin{equation}
		\label{generalj}
		\tilde{P}_{j, NIG}\stackrel{d}{=}V_j \prod_{i=1}^{j-1} (1-V_i)\nonumber
	\end{equation}
	for 
	\begin{equation}
		\label{pit_j}
		V_{j}|X_1, \dots, X_{j-1}=\frac{X_j}{T^{-1} +S_j}.
	\end{equation}
	
	\proof{  From Pitman (2003, eq. 97) and the result in Proposition 1
		\begin{equation}
			\label{stickPIT2}
			\tilde{P}_{j, NIG}\stackrel{d}{=} \frac{T^{-1}}{T^{-1}+S_{j-1}}- \frac{T^{-1}}{T^{-1}+S_j}.
		\end{equation}
		For $j=2$, then an easy manipulation of (\ref{stickPIT2}) yields
		\begin{equation}
			\label{p2}
			\tilde{P}_{2,NIG}\stackrel{d}{=}\frac{T^{-1}}{T^{-1}+X_1}\left(1 -\frac{T^{-1}+X_1}{T^{-1}+X_1+X_2}\right) \nonumber
		\end{equation}
		which can be rewritten as
		\begin{equation}
			\label{p2}
			\tilde{P}_{2, NIG}\stackrel{d}{=}\left(\frac{X_2}{T^{-1}+X_1+X_2}\right)\left(1- \frac{X_1}{T^{-1}+X_1}\right) \nonumber
		\end{equation}
		for $V_2=\left(\frac{X_2}{T^{-1}+X_1+X_2}\right)$ and $(1-V_1)=\left(1- \frac{X_1}{T^{-1}+X_1}\right)$.
		For $j> 2$ the result easily follows along the same lines. 
		\qed
	}
	
	{Notice that the stick-breaking construction obtained in Proposition 1 and Proposition 2 appears to provide both analytical and  computational advantages with respect to (\ref{struFAV}) and (\ref{vjFAV}), both to describe the sequence $(V_j)_{j \geq 1}$ of the random sticks and to model the dependence among them. Indeed, in place of a series of dependent generalized inverse Gaussian r.v.s plus a series of independent Stable ($1/2$) r.v.s, it just relies on a series of independent Chi-square r.v.s plus a single Inverse Gaussian r.v.. The possibility to exploit the {\it conditional} stick-breaking result in Aldous and Pitman (1998) comes from the representation of the NIG prior as a {\it mixed} Poisson-Kigman model.
		Conversely Favaro {\it et al.} (2012) move from the definition of the NIG model as a normalized subordinator, or, which is the same,  as an {\it unmixed} Poisson-Kingman model driven by the exponentially tilted $1/2$-stable L\'evy measure.
		As a consequence,  their stick-breaking representation is obtained resorting to Theorem 2.1 in Perman {\it et al.} (1992). 
		
		Note that the equivalence in distribution PK$(\rho_{\alpha}, \gamma^v)\stackrel{d}{=}$PK$(\rho_{\alpha}^{v})$, for $\rho_{\alpha}^v(s)=\rho_{\alpha}(s)\exp(-v^2s)$  the L\'evy density of an exponentially tilted $\alpha$-stable r.v., for $\alpha \in (0,1)$ is well-known to hold true (see Pitman, 2003, eq. 46, Sect. 4.2). 
	}
	
	
	To establish a link between Propositions 1 and 2 and the results in Favaro {\it et al.} (2012) we provide a description of the results in Proposition 2 in terms of Generalized Inverse Gaussian r.v.'s.  
	{\rema First recall that if $Y \sim IG (\mu, \lambda)$ then   $Y \sim GIG (\lambda/\mu^2, \lambda, -1/2)$ where for $x>0$ the density of $X \sim GIG(a, b, p)$ is given by
		\begin{equation}
			\label{gigdens}
			f_{GIG}(x)=\frac{(a/b)^{p/2}}{2K_p(\sqrt{ab})}x^{p-1}\exp(-(ax+b/x)/2)\nonumber,
		\end{equation}
		for $a>0$, $b>0$ and $p \in (-\infty, \infty)$. 
		From Proposition 1 $\tilde{P}_1\stackrel{d}{=}\frac{T}{T+1/X}$	for $X \sim \chi_1$ therefore $(1-\tilde{P}_1)^{-1}\stackrel{d}{=}(TX +1)$.
		Now recall that if $Y \sim IG(a, 1)$, for $a=E(Y)>0$ then, for $c>0$, $c Y \sim IG(ca, c)$ hence if $T \sim IG(1/v, 1)$ then $T(1-V_1)|V_1 \sim IG(\frac{1-v_1}{v}, 1-v_1)$. Thus 
		\begin{equation}
			T(1-V_1)|V_1 \sim GIG \left(\frac{v^2}{1-V_1}, 1-V_1, -1/2\right)\nonumber
		\end{equation}
		and the results in Propositions 1 and 2 can be restated in terms of Generalized Inverse Gaussian random variables as follows
		$$
		\begin{gathered}
			\label{miaGIG}
			V_1\stackrel{}{=}\frac{G_1}{G_1+Z_1} \quad \text { for } G_1 \sim \operatorname{GIG}\left(v^2, 1,-\frac{1}{2}\right), \\
			\left(V_i \mid V_1, \ldots, V_{i-1}\right)\stackrel{}{=}\frac{G_i}{G_i+Z_i} \quad \text { for } G_i \sim \operatorname{GIG}\left(\frac{v^2}{\prod_{j=1}^{i-1}(1-V_j)}, \prod_{j=1}^{i-1} (1-V_j),-\frac{1}{2}\right), \quad i \geq 2.
		\end{gathered}
		$$  
		and  $(Z_i)_i$ iid L\'evy $(0,1)$  r.v.s.   
		For example, for $i=2$, by (\ref{pit_j}) we get
		\begin{eqnarray}
			\label{v2}
			V_2|X_1 \stackrel{}{=}\frac{X_2}{\frac{1}{T}+X_1+X_2}\stackrel{}{=}\frac{T X_2}{TX_1 +1+TX_2}\nonumber
		\end{eqnarray}
		or equivalently 
		\begin{eqnarray}
			\label{v2}
			V_2|V_1 \stackrel{}{=}\frac{TX_2}{\frac{1}{(1-V_1)}+TX_2}\stackrel{}{=}\frac{T (1-V_1)}{T(1-V_1) +\frac{1}{X_2}}.\nonumber
		\end{eqnarray}
		Note that the discrepancy with respect to (\ref{vjFAV}) and (\ref{vjFAV2}), where
		\begin{equation}
			G_i \sim \operatorname{GIG}\left(\frac{v}{\prod_{j=1}^{i-1}\left(1-V_j\right)}, 1,-\frac{i}{2}\right)\nonumber,
		\end{equation} 
		arises since the treatment in Favaro et al. (2012) is in terms of the single parameter parametrization of the Inverse Gaussian r.v $X \sim IG (v)$ with $E(X)=v$, whose density and corresponding L\'evy density are given by  
		\begin{equation}
			f_{IG}(x)=\frac{v}{\sqrt{2\pi}} x^{-3/2}\exp\left\{-\frac{1}{2}\left(\frac{v^2}{x} +x\right)+v\right\}\nonumber
		\end{equation}
		and
		\begin{equation}
			\label{levyFAV}
			\rho_{IG}(s)=\frac{v}{2\pi} s^{-3/2}\exp(-s/2)\nonumber.
		\end{equation}
		Here instead, following Pitman (2003), we are dealing with $Y \sim IG(\mu, \lambda)$, $E(Y)=\mu=1/v$ and $\lambda=1$, with density  (\ref{IGdensity}), whose L\'evy density arises by tilting $\rho_{1/2}(s)=\frac{1}{2\sqrt{\pi}}s^{-3/2}$ by $\exp(-v^2 s)$. Overall this means that $X$ and $Y$ are linked by the transformation $X=Y/\mu_Y^2=v^2Y$.
	}

	
	A remarkable feature of the stick-breaking construction illustrated in Propositions 1 and 2 is that the mixing density $\gamma$ of the Poisson-Kingman model generating the inverse Gaussian model, namely the $IG(1/v, 1)$ density, appears explicitly in (\ref{struPIT}) and (\ref{pit_j}). This implies our results can be easily generalized to the size-biased permutation of the ranked atoms of {\it any} random discrete distribution in the class $PK(\rho_{1/2}, \gamma)$, that is, to the whole family of Gibbs-type priors of index $1/2$ (Gnedin and Pitman, 2005) and general mixing density $\gamma$. 
	
	A previous solution to this problem can be found in Favaro {\it et al.} (2014) that prove the following construction for the sticks of the $PK(1/2, \gamma_{T})$ family 
	\begin{equation}
		\left(
		V_{\frac12,i}
		\,\middle|\,
		V_{\frac12,1},\ldots,V_{\frac12,i-1},T_{}
		\right)
		\overset{}{=}
		\frac{X_i}{X_i+Y_i}\nonumber
	\end{equation}
	for
	\begin{equation}
		X_i^2 \sim {Ga}\!\left(\frac34,1\right) \nonumber
	\end{equation}
	and
	\begin{equation}
		Y_i^2 \sim IGa \left(
		\frac{1}{4},
		\frac{\frac{1}{
				4^{2}\, T_{}^{2}}}{
			\prod_{j=1}^{i-1}
			\left(1 - V_{\frac12,j}\right)^2
		}
		\right),\nonumber
	\end{equation}	
	with the convention that the empty product is defined to be unity and with
	${Ga}$ and $IGa$ denoting the Gamma distribution and the inverse Gamma distribution, respectively.
	Nevertheless a simpler description can be obtained relying on the results in Propositions 1 and 2.

	{\coro Let $(X_i)_{i \geq 1}$ be a sequence of iid $\chi_1^2$ r.v.s,  $S_j:=\sum_{i=1}^j X_i$, and $T$ a r.v. in $(0, \infty)$,  independent from $(X_i)_{i \geq 1}$, with density $\gamma_T(t)$, then the size-biased permutation $(\tilde{P}_j)$ of the ranked atoms $(P)^{\downarrow}_{j\geq 1} \sim PK(\rho_{1/2}, \gamma_T)$ admits the following stick-breaking construction 
		\begin{equation}
			\label{struPIT_gen}
			\tilde{P}_{1, T}\stackrel{d}{=}\frac{X_1} {T^{-1} + X_1}\nonumber
	\end{equation}}
	and
	\begin{equation}
		\label{generalj}
		\tilde{P}_{j, T}\stackrel{d}{=}V_j \prod_{i=1}^{j-1} (1-V_i)\nonumber
	\end{equation}
	for 
	\begin{equation}
		\label{pit_j_2}
		V_{j, T}|X_1, \dots, X_{j-1}, T\stackrel{}{=}\frac{X_j}{T^{-1} +S_j}\nonumber.
	\end{equation} 
	
	For a generalization of the stick-breaking construction of the NIG prior to the class of homogeneous normalized random measures see Favaro {\it et al.} (2016). Note that this generalization does not cover the general family $PK(\rho_{1/2}, \gamma_T)$, since the NIG prior is the only element in this class admitting a construction via normalization of a subordinator. 
	
	A general formula for the density of the stick-breaking weights of mixed Poisson-Kingman models driven by the stable $(\alpha)$ model for $\alpha \in (0,1)$ can be derived from formula (57) in Pitman (2003). See also Favaro and Walker (2013) for an application in the setting of slice sampling for mixture models. Note that the results in Propositions 1 and 2 cannot be extended to the general $\alpha$ case (see Miermont and Schweinsberg, 2003) since the analogy (\ref{analogy}) between the conditional Brownian excursion partition and the $PK(\rho_{1/2}|\alpha t)$ model is particular to the $\alpha=1/2$ case and does not extend to general $\alpha$.
	
	\subsection{James's (2013) stick-breaking construction}

	An additional easy-to-sample stick-breaking representation of the size-biased permutation of the ranked atoms of the normalized Inverse Gaussian random discrete distribution arises as a particular case of a general result derived in James (2013) for the family of {\it Poisson-Gamma $PG(\alpha, \zeta)$ partition models}. Those are random discrete distributions built upon a generalization of the construction of the $PD(\alpha, \theta)$ model as a random time changed Generalized Gamma subordinator, as in Pitman and Yor (1997, Prop. 21). See also James (2019).
	
	First recall that for $\alpha \in (0,1)$ and $\theta > -\alpha$ the {\it two parameter $(\alpha, \theta)$ Poisson-Dirichlet} distribution (Pitman, 1995, Pitman and Yor, 1997) is the law of the ranked atoms of a random discrete distribution $P$ which induces an exchangeable partition probability function (EPPF) in the form
	\begin{equation}
		\label{EPPFpd}
		p_{\alpha, \theta}(n_1, \dots, n_k)= \frac{(\theta +\alpha)_{k-1 \uparrow \alpha}}{(\theta +1)_{n-1}} \prod_{j=1}^k (1 -\alpha)_{n_j-1}\nonumber,
	\end{equation}
	for $(x)_{y \uparrow \alpha}=x(x+\alpha)\cdots(x+(y-1)\alpha)$ generalized rising factorials, $\sum_j n_j=n$ and $n_j$ for $j=1, \dots, k$ the sizes of the partition's blocks. The stick-breaking representation of the size-biased permutation of its ranked atoms  $P_{\alpha, \theta}=(\tilde{P}_{j}^{\alpha, \theta})$  is well-known to be $\tilde{P}_1^{\alpha, \theta}\stackrel{d}{=}V_1 \sim Beta(1-\alpha,\theta +\alpha)$ and $\tilde{P}_j^{\alpha, \theta}\stackrel{d}{=}V_j \prod_{j=1}^k (1-V_j)$ for independent $V_j \sim Beta(1-\alpha, \theta +j\alpha)$. 
	Now let $G_a$ denote a Gamma r.v. with shape parameter $a$ and scale $1$, $S_\alpha$ a positive Stable r.v. of index $\alpha \in (0,1)$ with density $f_\alpha(t)$ and $E_1$ an Exponential $(1)$ r.v.. Let $T_{\alpha,1}^{GG}(s): s \geq 0$ be a Generalized Gamma subordinator of parameters $(\alpha, 1)$ such that at time one $T_{\alpha, 1}^{GG}(1)$ has density $e^{-(t-1)}f_\alpha(t)$ then, by Proposition 21 in Pitman and Yor (1997), for $\theta>0$ and $J_1^{\downarrow} \geq J_2^{\downarrow} \geq J_3^{\downarrow} \dots$ its ranked jumps  over the random interval $[0, G_{\theta/\alpha}]$, then
	$
	T_{\alpha, 1}^{GG}(G_{\theta/\alpha})=\sum_{i=1}^\infty J_i \stackrel{d}{=} G_\theta
	$
	and is independent of the normalized ranked jumps $(J_i^{\downarrow}/T_{\alpha,1}^{GG}(G_{\theta/\alpha}))_i \sim PD(\alpha, \theta)$.
	
	
	Relying on this definition of the two-parameter Poisson-Dirichlet model $PD(\alpha, \theta)$, for $\theta >0$ as a {\it random} time changed subordinator of the generalized Gamma $(\alpha, 1)$ model James (2013) derived an alternative representation of its stick-breaking construction that is summarized in the following Proposition.

	{\prop (Corollary 8.2 in James, 2013) Let $(P)_j^{\downarrow} \sim PD(\alpha, \theta)$ then the random variables $W_j$ defining the stick-breaking representation of its size-biased permutation  $\tilde{P}_j=W_j\prod_{i=1}^{j-1}(1- W_i)$ admit the following alternative representation for $B_{1-\alpha, \alpha}^{(i)}$ IID Beta $(1-\alpha, \alpha)$, independent of $G_{\theta/\alpha+i-1}$ $\sim$ Gamma $(\theta/\alpha+i-1, 1)$, and $E_i$ IID exponential $(1)$ 
		\begin{equation}
			\label{stickLJ_pd}
			W_i^{\alpha, \theta}= B_{1-\alpha, \alpha}^{(i)}\left[1-\left(\frac{G_{\theta/\alpha} +G_{i-1}}{G_{\theta/\alpha}+G_{i-1}+E_i}\right)^{1/\alpha}\right].
		\end{equation}
	}
	
	{\it Proof}. This follows easily from elementary theory of Gamma and Beta r.v.s recalling that if $X_1 \sim Beta (a_1, b_1)$ and $X_2 \sim Beta (a_2, b_2)$ are independent and $a_2+b_2=a_1$ then $Y=X_1*X_2 \sim Beta (a_2, b_1+b_2)$. Additionally if $X_1 \sim Beta (a, 1)$ then $X_1^{1/\alpha} \sim Beta (a\alpha, 1)$. Therefore 
	$$
	W_i^{\alpha, \theta}= B_{1-\alpha, \alpha}^{(i)}\left[1-B_{\theta+(i-1)\alpha, 1}\right]=B_{1-\alpha, \theta+i\alpha}.
	$$	\qed
	
	Now, the following equivalence in distribution holds true when applying {\it non-random} spatial and temporal scaling to a generalized Gamma subordinator  $T_{\alpha, 1}^{GG}$
	\begin{equation}
		\label{invtempspat}
		\frac{T_{\alpha,1}^{GG} (v)}{v^{1/\alpha}}\stackrel{d}{=}T_{\alpha, v^{1/\alpha}}^{GG}(1)\nonumber
	\end{equation}
	with density $f_{\alpha, v^{1/\alpha}}^{GG}(t)= f_\alpha(t)\exp[-(v^{1/\alpha}t+v)]$ or alternatively in Pitman's (2003) notation for $v=b^{\alpha}$	
	with density $f_{\alpha, b}^{GG}(t)= f_\alpha(t)\exp[-(bt+b^{\alpha})]$.  
	
	It follows  that a normalized generalized Gamma subordinator $GG(\alpha, b)$ can be seen as a non-random {\it spatial} and {\it temporal} change of a normlized $GG (\alpha, 1)$ subordinator and a stick-breaking representation $\tilde{P}_{j}^{\alpha, b}=W_j\prod_{i=1}^{j-1}(1-W_j)$ of the size-biased permutation of its ranked atoms can be derived by (\ref{stickLJ_pd}) as having weights
	$$
	W_i^{\alpha, b}= B_{1-\alpha, \alpha}^{(i)}\left[1-\left(\frac{b^{\alpha} +G_{i-1, 1}}{b^{\alpha} +G_{i-1, 1}+E_i}\right)^{1/\alpha}\right].
	$$	
	An application of the result in Proposition 8.3 in James (2013) for general $\alpha \in (0,1)$ to the case $\alpha=1/2$  provides an additional constructive definition of the structural distribution of the normalized inverse Gaussian prior as summarized in the following Proposition.
	
	\prop  Let $B_{1/2, 1/2}$ be a Beta r.v. of parameters $(1/2, 1/2)$ and $G_v$ an independent r.v. in $(0,1)$ for which the following equality in distribution holds
	$$
	G_v\stackrel{d}{=}\left(1-\left(\frac{v}{v+E_1}\right)^2\right)
	$$
	for $E_1$ an Exponential $(1)$ and $v>0$. Then, according to James (2013), the distribution of the  first size-biased pick $\tilde{P}_{1, IG}$ from the ranked atoms of the normalized Inverse Gaussian partition model of parameter $
	(1/v, 1)$ admits the following representation  
	\begin{equation}
		\tilde{P}_{1, IG}\stackrel{d}{=}B_{1/2, 1/2}G_{v}
	\end{equation}
	whose density is given by
	\begin{eqnarray}
		\label{densLAN}
		f_{\tilde{P}_1}(p)=\frac{1}{2\pi} \frac{\exp(v)}{\sqrt{p(1-p)}} \int_p^1 \frac{\sqrt{1-p}}{\sqrt{g-p}}v(1-g)^{-3/2}\exp\{-v(1-g)^{-1/2}\}dg.
	\end{eqnarray}
	Additionally, for $\chi_1^2$ a Chi-squared r.v. and $T \sim IG(1/v, 1)$ then
	\begin{equation}
		\label{equiv}
		\frac{\chi_1^2}{\frac1T +\chi_1^2} \stackrel{d}{=} B_{1/2, 1/2}G_{v}
	\end{equation}
	
	\proof: By elementary probability the distribution function of $G_v$ corresponds to
	$$
	F_{G_v}(g)= 1-\exp\{-v[(1-g)^{-1/2} -1]\}
	$$
	with density 
	$$
	f_{G_v}(g)= \frac{v}{2}(1-g)^{-3/2} \exp\{-v[(1-g)^{-1/2}-1]\}dg.
	$$
	By scale mixture of $B_{1/2, 1/2}$ and $G_v$, for
	$
	f_{B_{1/2, 1/2}}(p)=\frac{1}{\pi}\frac{1}{\sqrt{p(1-p)}}
	$
	then
	$$
	f_{\tilde{P}_1}(p)=\int_p^1 f_{G_v}(y)f_{B_{1/2, 1/2}}(p/g)\frac{1}{g}dg
	$$
	hence
	$$
	f_{\tilde{P}_1}(p)=\int_p^1 \frac{1}{\pi}\frac{1}{\sqrt{(p/g)(1-p/g)}}\frac{1}{g}\frac{v}{2}(1-g)^{-3/2}\exp\{-v[(1-g)^{-1/2}-1]\}dg.
	$$
	Multiplying and dividing by $\sqrt{1-p}$ yields (\ref{densLAN})
	$$
	\label{lancdens}
	f_{\tilde{P}_1}(p)=\frac{1}{2\pi}\frac{v\sqrt{1-p} \exp{v}}{\sqrt{p(1-p)}}
	\int_p^1\frac{1}{\sqrt{g-p}}(1-g)^{-3/2}\exp\{-v[(1-g)^{-1/2}]\}dy.
	$$
	Now, to prove  (\ref{equiv})  it is enough to verify that
	$$
	\frac{v}{2} \int_p^1 \frac{\sqrt{1-p}}{\sqrt{g-p}}(1-g)^{-3/2} \exp \left\{-v(1-g)^{-1/2}\right\}dg=
	\int_0^\infty \frac{1}{2y^2}\exp\left\{-\frac{1}{2}\left(\frac{v^2 y}{1-p} +\frac{1}{y} \right)\right\}dy
	$$
	Let $1-g=z$  and $1-p=q$ then the left hand side corresponds to
	$$
	\frac{v}{2}\int_0^q \left(1-\frac{z}{q}\right)^{-1/2}  z^{-3/2}\exp\{-vz^{-1/2}\}dz
	$$
	then put $x=z/q$ to obtain
	$$
	\frac{v}{2}\int_0^1 (1-x)^{-1/2} (x)^{-3/2} q^{-1/2} \exp\{-v(qx)^{-1/2}\}dx
	$$
	and again put $y=1/x -1$, $x=1/(1+y)$ then
	$$
	\frac{v}{2}\int_0^\infty y^{-1/2} q^{-1/2}\exp\left\{-v(q/(1+y))^{-1/2}\right\}dy=\frac{v}{\sqrt{q}}\int_0^\infty \frac{1}{2\sqrt{y}} \exp\{-v \left(\frac{q}{1+y}\right)^{-1/2}\}dq=
	$$
	$$
	=\frac{v}{\sqrt{q}} K_1(\frac{v}{\sqrt{q}}).
	$$
	Now, it can be numerically shown that
	\begin{equation}
		\label{besselINTE}
		\int_0^\infty \frac{1}{2\sqrt{y}} \exp\left\{-v \left(\frac{q}{1+y}\right)^{-1/2}\right\}dq= K_1\left(\frac{v}{\sqrt{q}}\right),\nonumber
	\end{equation}  
	by (\ref{pit_bes}) and (\ref{pit_ig}) the proof is completed. \qed

	Again from Proposition  8.3 in James (2013) the result for $\tilde{P}_{1, IG}$ in the previous Proposition extends to the whole sequence of the stick-breaking weights as follows.

	{\prop  [James, 2013, Th. 8.1] For $i=1, \dots, j$ let $(B_{1/2, 1/2}^{(i)})_i$ be independent Beta r.v.'s and  $G_{i, 1}=\sum_{l=1}^i E_l$ for $(E_l)_l$ independent exponential $(1)$ r.v.'s, independent from the  $(B_{1/2, 1/2}^{(i)})_i$, then the following stick-breaking construction for the size-biased permutation of the ranked atoms of the normalized Inverse Gaussian prior holds true
		\begin{equation}
			\tilde{P}_{j, NIG}\stackrel{d}{=}V_j \prod_{i=1}^{j-1} (1-V_i)\nonumber
		\end{equation}
		for 
		\begin{equation}
			\label{generalj2}
			V_j^{NIG}= B_{1/2, 1/2}^{(j)}\left[1-\left(\frac{v +G_{j-1, 1}}{v +G_{j-1, 1}+E_j}\right)^{2}\right]
		\end{equation}
	}

	A detailed proof for (\ref{generalj2}) and (\ref{pit_j}) to be equal in distribution is beyond the scope of this paper.  A simulation-based check can be easily done.
	
	\medskip
	
	The results in Propositions 4 and 5 can be easily generalized to a larger family of random discrete distributions, the one generated by applying a  {\it random} spatial and temporal change to the Inverse Gaussian subordinator, i.e. randomizing the parameter $v$ in (\ref{generalj2}).  This operation may disclose the possibility to investigate a whole new class of Bayesian nonparametric priors, since the availability of an easy-to-sample stick-breaking construction of the size-biased permutation of its ranked atoms would address the computational issues that usually prevent the implementation of conditional and marginal simulation algorithms for posterior inference. This line of research will be the topic of some future papers.
	
	\section{References}

	\begin{enumerate} 
		\item Aldous, D. Pitman, J. (1998) The standard additive coalescent. {\it Ann. Probab.} 26 (4) 1703-1726.
		\item Canale, A., Corradin, R., Nipoti, B. (2022). Importance conditional sampling for Pitman Yor mixtures. {\it Statistics and Computing}, 32(3) 1-18
		\item Favaro, S. Walker, S.G. (2013) Slice Sampling $\sigma$-Stable Poisson-Kingman Mixture Models. {\it JCGS}, 22, 4, 830-847
		\item Favaro, S., Lijoi, A., and Pr\"unster, I. (2012) On the Stick-Breaking Representation of Normalized Inverse Gaussian Priors.  {\it Biometrika}, 99, 663-674.
		
		\item Favaro, S., Lomeli, M., Nipoti, B. and Teh, Y. W. (2014) On the Stick-breaking representation of $\sigma$-stable Poisson-Kingman models. {\it Electron. J. Statist.} 8 (1) 1063 - 1085.
		\item Favaro, S. Lijoi, A. Nava, C. Nipoti, B. Pr\"unster, Teh, Y.I.W. (2016) On the Stick-Breaking Representation for Homogeneous NRMIs. {\it  Bayesian Analysis} 11 (3) 697 - 724
		\item Ferguson, T.S. (1973) A Bayesian Analysis of Some Nonparametric Problems. {\it Ann. Statist.} 1 (2) 209 - 230. https://doi.org/10.1214/aos/1176342360
		
		\item Gnedin, A., and Pitman, J. (2005), Exchangeable Gibbs Partitions and Stirling Triangles, Zapiski Nauchnych Seminarov POMI, 325, 83-102. [834]
		\item James, L. F. (2013) Stick-breaking PG($\alpha, \zeta$)-Generalized Gamma Processes. {\it arXiv:1308.6570} [math-PR]
		\item James, L. F. (2019) Stick-breaking Pitman-Yor processes given the species sampling size. {\it arXiv:1908.07186 [math.ST] }
		\item Lijoi, A., Mena, R. H., and Pr\"unster, I. (2005) Hierarchical Mixture Modelling With Normalized Inverse-Gaussian Priors. {\it JASA}, 100, 1278-1291. 
		\item Miermont, G., and Schweinsberg, J. (2003) Self-similar fragmentations and stable subordinators. In Sèminaire de Probabilitè XXXVII, vol. 1832 of LNh in Maths. Springer, Berlin, pp. 333-359.
		\item Perman, M., Pitman, J. and Yor, M. (1992) Size-biased sampling of Poisson point processes and excursions. {\it  Probab. Th. Rel. Fields} 92, 21-39. 
		\item Pitman, J. (1995) Exchangeable and partially exchangeable random partitions. {\it  Probab. Th. Rel. Fields} 102, 145-158.
		
		\item  Pitman, J. (2003) Poisson-Kingman Partitions, in Science and Statistics: A Festschrift for Terry Speed, ed. D. R.
		Goldstein, Beachwood, OH: IMS, pp. 1-34.
		\item Pitman and M. Yor. (1992) Arcsine laws and interval partitions derived from a stable subordinator. Proceedings of the London Mathematical Society (3), 65:326-356, 1992. 
		\item Pitman, J., and Yor, M. (1997), The Two Parameter Poisson-Dirichlet Distribution Derived From a Stable
		Subordinator. {\it Annals of Probability}, 27, 1870-1902			
	\end{enumerate}
	
\end{document}